\documentclass{amsart}

\usepackage{amsmath,amsfonts}

\usepackage{amscd,amsthm}

\theoremstyle{plain} 

\newtheorem{theorem}{Theorem}[section]

\newtheorem{lemma}[theorem]{Lemma}

\newtheorem{proposition}[theorem]{Proposition}

\newtheorem{definition}[theorem]{Definition}

\newtheorem{Proposition 3.9}{prop39}

\def\limind{\mathop{\oalign{lim\cr\hidewidth$\longrightarrow$\hidewidth\cr}}}

\numberwithin{equation}{section}

\begin{document}

\vspace{-30em}
IN PRESS: BULLETIN SOCIETE MATHEMATIQUE DE FRANCE
\vspace{+3em}

\title[Linearization of germs of diffeomorphisms]
{Linearization of analytic and non--analytic germs of 
diffeomorphisms of $({\mathbb C},0)$}

\author{Timoteo Carletti, Stefano Marmi}

\date{\today}

\address[Timoteo Carletti, Stefano Marmi]{Dipartimento di Matematica 
``Ulisse Dini'',  Viale 
Morgagni 67/A, 50134 Florence, Italy}

\email[Timoteo Carletti]{carletti@udini.math.unifi.it}
\email[Stefano Marmi]{marmi@udini.math.unifi.it}


\keywords{Siegel's center problem, 
small divisors, Gevrey classes}

\begin{abstract}
We study 
Siegel's center problem on the linearization of germs of 
diffeomorphisms in one variable. In addition of the classical 
problems of formal and analytic 
linearization, we give sufficient conditions for the linearization to 
belong to some algebras of ultradifferentiable 
germs closed under composition and derivation, including Gevrey classes.  

In the analytic case we give a positive answer to a question of J.-C. 
Yoccoz on the optimality of the estimates obtained by 
the classical majorant series method. 

In the ultradifferentiable  case we prove  
that the Brjuno condition is sufficient for the linearization 
to belong to the same class of the germ. If one allows the 
linearization to be less regular than the germ one finds new 
arithmetical conditions, weaker than the Brjuno condition. 
We briefly discuss the optimality of our results. 
\end{abstract}

\maketitle

\section{introduction}

\indent 
In this paper we study the Siegel center problem
\cite{Herman}. 
Consider two subalgebras $A_{1}\subset A_{2}$ of 
$z\mathbb{C}\left[\left[ z\right]\right]$
closed with respect to the composition of formal series. For example 
$z\mathbb{C}\left[\left[ z\right]\right]$, 
$z\mathbb{C}\{  z \}$ (the usual analytic case) or 
Gevrey--$s$ classes, $s>0$ (i.e. series $F(z)=\sum_{n\ge 0}f_{n}z^n$
such that there exist $c_{1}, c_{2}>0$ such that 
$|f_{n}|\le c_{1}c_{2}^n(n!)^s$ for all $n\ge 0$). Let 
$F\in A_{1}$ being such that $F^{\prime}\left( 0 \right)=\lambda\in 
\mathbb{C}^{*}$. We say that $F$ is linearizable in $A_{2}$ if there 
exists $H\in A_{2}$ tangent to the identity and such that 
\begin{equation}
F \circ H  = H \circ R_\lambda  
\label{homological}
\end{equation}
where $R_\lambda \left( z \right)= \lambda z$.
When $|\lambda|\not= 1$, the Poincar\'e-Konigs linearization theorem 
assures that $F$ is linearizable in $A_{2}$. When $|\lambda|=1$, 
$\lambda = e^{2\pi i \omega}$, the 
problem is much more difficult, especially if one looks for 
{\it necessary and sufficient} conditions on $\lambda$ which assure that {\it all} 
$F\in A_{1}$ {\it with the same $\lambda$ are linearizable in 
$A_{2}$}. The only trivial case is $A_{2}= z\mathbb{C}\left[\left[ 
z\right]\right]$ (formal linearization) for which one only needs to 
assume that $\lambda$ is not a root of unity, i.e. 
$\omega\in \mathbb{R}\setminus\mathbb{Q}$.

In the analytic case $A_{1}=A_{2}=z\mathbb{C}\{z\}$ let 
$S_{\lambda}$ denote the space of analytic germs $F\in z\mathbb{C}\{z\}$
analytic and injective in the unit disk $\mathbb{D}$ and such that 
$DF(0)=\lambda$ (note that any $F\in z\mathbb{C}\{z\}$ tangent to 
$R_{\lambda}$ may be assumed to belong to $S_{\lambda}$ provided that 
the variable $z$ is suitably rescaled). Let 
$R(F)$ denote 
the radius of convergence of the unique tangent to the identity linearization $H$ associated 
to $F$.  
J.-C. Yoccoz \cite{Yoccoz} proved that the {\it Brjuno condition} 
(see Appendix \ref{appfrazioni}) is necessary and sufficient for having $R(F)>0$ 
for all $F\in S_{\lambda}$. More precisely  
Yoccoz proved the following estimate:  assume that $\lambda = e^{2\pi 
i \omega}$ is a 
Brjuno number. There exists a universal constant $C>0$ (independent 
of $\lambda$) such that 
\begin{equation*}
\label{Yoccozestimate}
|\log R(\omega ) +B(\omega )|\le C
\end{equation*}
where  $R(\omega )= \inf_{F\in S_{\lambda}} R(F)$ and $B$ is the 
Brjuno function (\ref{Brjunofunction}).  Thus 
$\log R(\omega )\ge  -B(\omega ) - C$.

\indent
Brjuno's proof \cite{Brjuno} gives an estimate of the form 
\begin{equation*}
	\label{Brjunoestimate}
	\log r(\omega )\ge -C'B(\omega )-C''
\end{equation*}
where one can choose $C'=2$ \cite{Herman}. 
Yoccoz's proof is based on a geometric renormalization argument and 
Yoccoz himself asked whether or not was possible to obtain $C'=1$ by 
direct manipulation of the power series expansion of the 
linearization $H$ as in Brjuno's proof (\cite{Yoccoz}, Remarque 
2.7.1, p. 21). 
Using an arithmetical lemma due to Davie
\cite{Davie} (Appendix \ref{appDavie}) we give a positive answer  
(Theorem \ref{Yoccozlower})
to Yoccoz's question.

\indent
We then consider 
the more general ultradifferentiable case $A_{1}\subset A_{2}\not= z{\mathbb C}\{z\}$. If one requires $A_{2}=A_{1}$, i.e.  the 
linearization $H$ to be as regular as the given germ $F$, once again the 
Brjuno condition is sufficient. Our methods do not allow us to 
conclude that the Brjuno condition is also necessary, 
a statement which is in general false as we show in section 2.3 where 
we exhibit a Gevrey--like class for which the sufficient condition 
coincides with the optimal arithmetical condition for the associated 
linear problem. Nevertheless it is 
quite interesting to notice that given any algebra of formal power 
series which is closed under composition (as it should if one whishes 
to study conjugacy problems) and derivation a germ in the 
algebra is linearizable {\it in the same algebra} if the Brjuno 
condition is satisfied. 

If the linearization is allowed to be less regular than the given germ 
(i.e. $A_{1}$ is a proper subset of $A_{2}$) 
one finds a new arithmetical condition, weaker than the Brjuno 
condition. This condition is also optimal if the small divisors are 
replaced with their absolute values as we show in section 2.4.
We discuss two examples, including Gevrey--$s$ classes.\footnote{We refer the reader interested in small divisors and Gevrey--$s$ classes to \cite{Lochak, GramchevYoshino1, GramchevYoshino2}.}

\indent
{\it Acknwoledgements.} 
We are  grateful to J.--C. Yoccoz for 
a very stimulating discussion concerning Gevrey 
classes and small divisor problems. 

\section{the Siegel center problem}
Our first step 
will be the formal solution of equation (\ref{homological}) assuming 
only that $F\in z\mathbb{C}[[z]]$. Since  $F\in 
z\mathbb{C}\left[ \left[ z \right] \right]$ is assumed to be  tangent to 
$R_{\lambda}$ then $F(z)= \sum_{n\ge 
1}f_{n}z^n$ with $f_{1}=\lambda$.  
Analogously since  $H\in z\mathbb{C}\left[ \left[ z \right] \right]$ 
is tangent to the identity $H(z)=\sum_{n=1}^\infty
h_{n}z^n$ with $h_{1}=1$. If 
$\lambda$ is not a root of 
unity equation (\ref{homological}) has a unique solution 
$H\in z^{}\mathbb{C}\left[ \left[ z \right] \right]$ tangent to the 
identity: the 
power series coefficients satisfy  the recurrence relation 
\begin{equation}
\label{recursive}
h_{1}= 1 \; , \; h_{n} = \frac {1}{\lambda^n -\lambda}
\sum_{m=2}^n f_{m}\sum_{n_{1}+\ldots +n_{m}= n\, , \, 
n_{i}\ge 1 } h_{n_{1}}\ldots h_{n_{m}}\; . 
\end{equation}

In \cite{Teo} it is shown how to generalize the classical Lagrange inversion 
formula to non--analytic inversion problems on the field of formal power series
so as to obtain an explicit non--recursive formula for the power series
coefficients of $H$.

\subsection{The analytic case: a direct proof of Yoccoz's lower bound}

Let $S_{\lambda}$ denote the space of germs $F\in z\mathbb{C}\{z\}$ analytic and 
injective in the unit disk $\mathbb{D}=\{z\in \mathbb{C}\, , \, 
|z|<1\}$ such that $DF(0)=\lambda$ and assume
that $\lambda = e^{2\pi i\omega}$ with $\omega\in 
\mathbb{R}\setminus\mathbb{Q}$. With the topology of uniform 
convergence on compact subsets of $\mathbb{D}$, $S_{\lambda}$ is 
a compact space. Let $H_{F}\in z\mathbb{C}[[z]]$ denote the unique tangent to the 
identity formal linearization associated to $F$, i.e. the unique 
formal solution of (\ref{homological}). Its power series coefficients 
are given by (\ref{recursive}). Let $R(F)$ denote the radius of 
convergence of $H_{F}$. Following Yoccoz (\cite{Yoccoz}, 
p. 20) we define 
\begin{equation*}
R(\omega ) = \inf_{F\in S_{\lambda}} R(F)\; .
\end{equation*}
We will prove the following 

\begin{theorem} {\bf Yoccoz's lower bound}.  
\label{Yoccozlower}
\begin{equation}
\log R(\omega) \ge -B(\omega ) -C
\end{equation}
where $C$ is a universal constant (independent of $\omega$) and $B$ is 
the Brjuno function (\ref{Brjunofunction}). 
\end{theorem}

Our method of proof of Theorem \ref{Yoccozlower} will be to apply  an arithmetical 
lemma due to Davie (see Appendix B) to estimate the small divisors 
contribution to (\ref{recursive}). This is actually a variation of 
the classical majorant series method as used in \cite{Siegel, Brjuno}. 

\proof 

Let $s\left( z\right) =\sum_{n\geq 1}s_nz^n$ be  the 
unique  solution analytic at $z=0$ of the equation 
$s\left( z\right) =z+\sigma \left( s\left( z\right) \right)$, where 
$\sigma (z) = \frac{z^{2}(2-z)}{(1-z)^{2}}= \sum_{n\ge 2}nz^n$. 
The coefficients satisfy 
\begin{equation}
\label{recursive-s}
s_{1}= 1 \; , \; s_{n} = 
\sum_{m=2}^n m\sum_{n_{1}+\ldots +n_{m}= n\, , \, 
n_{i}\ge 1 } s_{n_{1}}\ldots s_{n_{m}}\; . 
\end{equation}
Clearly there exist two positive constants $\gamma_{1},\gamma_{2}$
such that 
\begin{equation}
\label{sestimate}
|s_n| \leq \gamma_{1}\gamma_{2}^{n}\; . 
\end{equation}
From the recurrence relation (\ref{recursive}) and 
Bieberbach--De Branges's bound
$|f_{n}|\le n $ for all $n\ge 2$ 
we obtain 
\begin{equation}
\label{hest1}
|h_{n}| \le \frac {1}{|\lambda^n -\lambda|}
\sum_{m=2}^n m\sum_{n_{1}+\ldots +n_{m}= n\, , \, 
n_{i}\ge 1 } |h_{n_{1}}|\ldots |h_{n_{m}}|\; . 
\end{equation}
We now deduce by induction on $n$ that $|h_{n}|\le s_{n}e^{K(n-1)}$
for $n\ge 1$, where $K$ is defined in Appendix B. If we assume this 
holds for all $n'<n$ then the above inequality gives 
\begin{equation}
\label{hest2}
|h_{n}| \le \frac {1}{|\lambda^n -\lambda|}
\sum_{m=2}^n m\sum_{n_{1}+\ldots +n_{m}= n\, , \, 
n_{i}\ge 1 } s_{n_{1}}\ldots s_{n_{m}}e^{K(n_{1}-1)+\ldots K(n_{m}-1)}
\; .
\end{equation}
But $K(n_{1}-1)+\ldots K(n_{m}-1)\le K(n-2)\le K(n-1)+\log |\lambda^n -\lambda|$
and we deduce that 
\begin{equation}
	\label{hest3}
	|h_{n}|\le e^{K(n-1)}\sum_{m=2}^n m\sum_{n_{1}+\ldots +n_{m}= n\, , \, 
n_{i}\ge 1 } s_{n_{1}}\ldots s_{n_{m}}=s_{n} e^{K(n-1)}\; , 
\end{equation}
as required. Theorem 2.1 then follows from the fact that 
$n^{-1}K(n)\le B(\omega )+\gamma_{3}$ for some universal constant 
$\gamma_{3}>0$ (Davie's lemma, Appendix B). 

\endproof

\subsection{The ultradifferentiable case}

A classical result of Borel says that the map $J_{\mathbb{R}}\, : \, 
\mathcal{C}^\infty ([-1,1],\mathbb{R}) \rightarrow \mathbb{R}[[x]]$ which 
associates to $f$ its Taylor series at $0$ is surjective. 
On the other hand, $\mathbb{C}\{z\} = 
\limind_{r>0}\mathcal{O}(\mathbb{D}_{r})$, where 
$\mathbb{D}_{r}=\{z\in \mathbb{C}\, , \, |z|<r\}$ and 
$\mathcal{O}(\mathbb{D}_{r})$ is the $\mathbb{C}$--vector space of 
$\mathbb{C}$--valued functions analytic in $\mathbb{D}_{r}$. 
Between $\mathbb{C}[[z]]$ and $\mathbb{C}\{z\}$ one has many 
important algebras of ``ultradifferentiable'' power series
(i.e. asymptotic expansions at $z=0$ of functions which are 
``between'' $\mathcal{C}^\infty$ and $\mathbb{C}\{z\}$). 

\indent 
In this part we will study the case $A_{1}$ or $A_{2}$ (or both) is neither 
$z\mathbb{C}\{  z \}$ nor $z\mathbb{C}\left[ \left[ z \right] 
\right]$ but a general ultradifferentiable algebra $z
\mathbb{C}\left[ \left[ z \right] \right]_{(M_{n})}$
defined as follows.

Let $(M_{n})_{n\ge 1}$ be a sequence of positive real numbers such 
that: 
\item{0.} $\inf_{n\ge 1} M_{n}^{1/n}>0$; 
\item{1.} There exists $C_{1}>0$ such that $M_{n+1}\le C_{1}^{n+1}M_{n}$ 
for all $n\ge 1$; 
\item{2.} The sequence $(M_{n})_{n\ge 1}$ is logarithmically convex; 
\item{3.} $M_{n}M_{m}\le M_{m+n-1}$ for all $m,n\ge 1$.

\begin{definition}
\label{ultradifferentiable}

Let $f= \sum_{n\ge 1}f_{n}z^n\in z\mathbb{C}
\left[ \left[ z \right] \right]$; $f$ belongs to 
the algebra 
$z\mathbb{C}\left[ \left[ z \right] \right]_{(M_{n})}$
if there exist two positive constants $c_{1},c_{2}$ 
such that 
\begin{equation}
\label{coeffs}
|f_{n}| \le c_{1}c_{2}^nM_{n}\;\; 
\hbox{for all}\; n\ge 1\; . 
\end{equation}
\end{definition}

The role of the above assumptions on the sequence 
$(M_{n})_{n\ge 1}$ is the following: 0. assures that $z\mathbb{C}\{  z 
\}\subset z\mathbb{C}\left[ \left[ z \right] \right]_{(M_{n})}$; 
1. implies that $z\mathbb{C}\left[ \left[ z \right] \right]_{(M_{n})}$
is stable for derivation. Condition 2. means that $\log M_{n}$ is 
convex, i.e. that the sequence $(M_{n+1}/M_{n})$ is increasing; it 
implies that $z\mathbb{C}\left[ \left[ z \right] \right]_{(M_{n})_{n\ge 1}}$
is an algebra, i.e. stable by multiplication. Condition 3. implies 
that this algebra is {\it closed for composition}: if $f,g\in 
z\mathbb{C}\left[ \left[ z \right] \right]_{(M_{n})_{n\ge 1}}$ then 
$f\circ g \in z\mathbb{C}\left[ \left[ z \right] \right]_{(M_{n})_{n\ge 
1}}$.  This is a very natural assumption since we will study a {\it 
conjugacy} problem. 

Let $s>0$. A very important example of ultradifferentiable algebra is given by 
the algebra of {\it Gevrey}--$s$ series which is obtained chosing 
$M_{n}= (n!)^s$. It is easy to check that the assumptions 0.--3. are 
verified. But also more rapidly growing sequences may be considered 
such as $M_{n} = n^{a n^b}$ with $a >0$ and 
$1<b <2$. 

We then have the following 

\begin{theorem}
\label{theultradifferentiable}

\item{1.} If $F\in z\mathbb{C}\left[ \left[ z \right] \right]_{(M_{n})}$
and $\omega$ is a Brjuno number then also the linearization $H$ 
belongs to the same algebra $z\mathbb{C}\left[ \left[ z \right] 
\right]_{(M_{n})}$.
\item{2.} If $F\in z\mathbb{C}\left\{ z\right\}$ and $\omega$ verifies 
\begin{equation}
\label{BrjunoM}
	\limsup_{n\rightarrow +\infty} \left( \sum_{k=0}^{k(n)}
	\frac{\log q_{k+1}}{q_{k}} - \frac{1}{n}\log M_{n}\right) <
	+\infty
\end{equation}
where $k(n)$ is defined by the condition $q_{k(n)}\le n< q_{k(n)+1}$, 
then the linearization $H\in z\mathbb{C}\left[ \left[ z \right] \right]_{(M_{n})}$.
\item{3.} Let  $F\in z\mathbb{C}\left[ \left[ z \right] 
\right]_{(N_{n})}$, 
where the sequence $(N_{n})$ verifies 0,1,2,3 and is asymptotically 
bounded by the sequence $(M_{n})$ (i.e. $M_{n}\ge N_{n}$ for all 
sufficiently large $n$). 
If $\omega$ verifies 
\begin{equation}
\label{BrjunoMN}
	\limsup_{n\rightarrow +\infty} \left( \sum_{k=0}^{k(n)}
	\frac{\log q_{k+1}}{q_{k}} - \frac{1}{n}\log \frac{M_{n}}{N_{n}}\right) <
	+\infty
\end{equation}
where $k(n)$ is defined by the condition $q_{k(n)}\le n< q_{k(n)+1}$, 
then the linearization $H\in z\mathbb{C}\left[ \left[ z \right] 
\right]_{(M_{n})}$.

\end{theorem}

Note that conditions (\ref{BrjunoM})  and (\ref{BrjunoMN}) 
are generally weaker than the 
Brjuno condition. For example if given $F$ analytic one only requires 
the linearization $H$ to be Gevrey--$s$ then one can allow the 
denominators $q_{k}$ of the continued fraction expansion of $\omega$ 
to verify $q_{k+1} = {\mathcal O} (e^{\sigma q_{k}}) $ 
for all $0< \sigma \le s$ whereas an exponential growth rate of the 
denominators of the convergents is clearly 
forbidden from the Brjuno condition. If the linearization is required 
only to belong to the class $z\mathbb{C}\left[ \left[ z \right] \right]_{(M_{n})}$
with $M_{n} = n^{a n^b}$,  with $a >0$ and 
$1<b <2$, one can even have $q_{k+1} = {\mathcal O} (e^{\alpha q_{k}^\beta}) $ 
for all $\alpha >0 $ and $1 <\beta <b$ and the series $\sum_{k\ge 0} \frac{\log 
q_{k+1}}{q_{k}^b}$ converges. This kind of series  have been studied in detail 
in \cite{MMY}. 

\proof
We only prove (\ref{BrjunoMN}) which clearly implies (\ref{BrjunoM}) 
(choosing $N_{n}\equiv 1$) and also assertion 1. (choosing 
$M_{n}\equiv N_{n}$). 

Since it is not restrictive to assume $c_{1}\ge 1$ and $c_{2}\ge 1$ in 
$|f_{n}|\le c_{1}c_{2}^nN_{n}$ one can immediately check by induction 
on $n$ that $|h_{n}|\le c_{1}^{n-1}c_{2}^{2n-2}s_{n}N_{n}e^{K(n-1)}$, 
where $s_{n}$ is defined in (\ref{recursive-s}). Thus by  
(\ref{sestimate}) and Davie's lemma one has 
\begin{equation*}
\frac{1}{n}\log\frac{|h_{n}|}{M_{n}}\le 
c_{3}+ \frac{1}{n}\log\frac{N_{n}}{M_{n}} + \sum_{k=0}^{k(n)}
	\frac{\log q_{k+1}}{q_{k}}
\end{equation*}
for some suitable constant $c_{3}>0$. 
\endproof

{\it Problem.} Are the arithmetical conditions stated in Theorem 
\ref{theultradifferentiable} optimal? In particular is it true that 
given any algebra $A= z\mathbb{C}\left[ \left[ z \right] 
\right]_{(M_{n})}$ and $F\in A$ then $H\in A$ {\it if and only if} 
$\omega$ is a Brjuno number?

We believe that this problem deserves further investigations and that 
some surprising results may be found. In the next two sections we will 
give some preliminary results.

\subsection{A Gevrey--like class where the linear and non 
linear problem have the same sufficient arithemtical condition}

Let $\mathbb{C} \left[ \left[z \right] 
\right]_s$ denote the algebra of Gevrey--$s$ complex formal power 
series, $s>0$. If $s^{\prime}>s>0$ then 
$z\mathbb{C} \left[ \left[z \right] \right]_s \subset z\mathbb{C} 
\left[ \left[z \right] \right]_{s^\prime}$; let 
\begin{equation*}
A_{s} = \bigcap_{s^\prime > s} z\mathbb{C} \left[ \left[z \right] 
\right]_{s^\prime}\; . 
\end{equation*}
Clearly $A_{s}$ is an algebra stable w.r.t. derivative and composition. 
This algebra can be equivalently characterized requiring that 
given $f \left( z \right) = \sum_{n \geq 1} f_n z^n
\in z\mathbb{C} \left[ \left[z \right] 
\right]$ one has 
\begin{equation}
\label{conditionAs}
\limsup_{n \rightarrow \infty} \frac{\log \lvert f_n \rvert}{n \log n} 
\leq s
\end{equation}

Consider Euler's derivative (see \cite{Du}, section 4)
\begin{equation}
(\delta_{\lambda}f)(z) = \sum_{n=2}^\infty (\lambda^n-\lambda)f_{n}z^n
\; , 
\end{equation}
with 
$\lambda = e^{2\pi i \omega}$. It acts linearly on $zA_{s}$ and it is 
a linear automorphism of $zA_{s}$ if and only if 
\begin{equation}
\label{lineararithmetic}
\lim_{k 
\rightarrow \infty}\frac{\log q_{k+1}}{q_k \log q_k}=0
\end{equation}
where, as usual, $\left( q_k \right)_{k \in \mathbb{N}}$ is the 
sequence of the denominators of the 
convergents of $\omega$. This fact can be easily checked by applying 
the law of the best approximation 
(Lemma \ref{bestapproximation}, Appendix ~\ref{appfrazioni})
and the charaterization (\ref{conditionAs}) to 
\begin{equation*}
h(z) = (\delta_{\lambda}^{-1}f)(z) = \sum_{n \geq 2} 
\frac{f_n}{\lambda^n -\lambda}z^n\; . 
\end{equation*}
Note that the arithmetical condition $\log q_{k+1}={\it o} \left(q_k \log q_k 
\right)$ is much weaker than Brjuno's condition.

We now consider the Siegel problem associated to a germ $F\in A_{s}$. 
Applying the third statement of Theorem \ref{theultradifferentiable}
with $N_n=\left( n! \right)^{s+\eta}$ 
and $M_n=\left( n! \right)^{s+\epsilon}$ for any positive fixed 
$\epsilon > \eta>0$ one finds that if the following arithmetical 
condition is satisfied
\begin{equation}
\label{nonlineararithmetic}
\lim_{k \rightarrow \infty}\frac{1}{\log q_{k}} \sum_{i=0}^{k}
\frac{\log q_{i+1}}{q_i}=0
\end{equation}
then the linearization $H_{F}$ also belongs to 
$A_{s}$.\footnote{In Theorem \ref{theultradifferentiable} we proved 
that a sufficient condition with this choice of $M_n$ and $N_n$ is 
\begin{equation*}
\limsup_{n\rightarrow +\infty} \left( \sum_{i=0}^{k(n)}
	\frac{\log q_{i+1}}{q_{i}} - \frac{\epsilon - \eta}{n} \log \left( n! 
\right) \right) \leq C < +\infty
\end{equation*}
which can be rewritten as
\begin{equation*}
\limsup_{n\rightarrow +\infty} \left( \sum_{i=0}^{k(n)}
	\frac{\log q_{i+1}}{q_{i}} - \left( \epsilon - \eta \right) \log 
q_{k\left( n \right)} -C \right) =0
\end{equation*}
from which 
(\ref{nonlineararithmetic}) is just obtained dividing by
$\log q_{k\left( n \right)}$.}

The equivalence of (\ref{nonlineararithmetic}) and 
(\ref{lineararithmetic}) is the object of the following 

\begin{lemma}
\label{equivalency}
Let $\left( q_l \right)_{l \geq 0}$ be the sequence of  
denominators of the convergents 
of $\omega \in \mathbb{R} \setminus \mathbb{Q}$. The following statements 
are all equivalent:
\begin{enumerate}
\item $\lim_{n \rightarrow \infty} \frac{1}{\log n}\sum_{l=0}^{k\left( 
n \right)}\frac{\log q_{l+1}}{q_l}=0$ 
\item $\sum_{l=0}^{k\left( n \right)}\frac{\log q_{l+1}}{q_l}={\it 
o}\left( \log q_{k} \right)$
\item $\log q_{k+1}={\it o}\left( q_k \log q_k \right)$
\end{enumerate}
\end{lemma}

\proof 

1. $\Longrightarrow$ 2. is trivial (choose $n=q_{k\left( n \right)}$). 

\indent

2. $\Longrightarrow$ 3.
Writing for short $k$ istead of $k \left( n \right)$
\begin{align*}
\frac{1}{\log q_k}\sum_{l=0}^{k} \frac{\log q_{l+1}}{q_l} &= \frac{\log 
q_{k+1}}{q_k \log q_{k}}+\frac{1}{\log q_k}\sum_{l=0}^{k-1} \frac{\log 
q_{l+1}}{q_l} \\
                &=\frac{\log q_{k+1}}{q_k \log q_{k}}+ \frac{{\it o}\left( 
\log q_{k-1} \right)}{\log q_{k}}
\end{align*}
Since $\lim_{k\rightarrow \infty} \frac{{\it o}\left( 
\log q_{k-1} \right)}{\log q_{k}} = 0$  we get 3.

\indent 3. $\Longrightarrow$ 1.
First of all note that since $q_{k\left( n 
\right)} \leq n$ 2.  trivially implies 1. Thus it is enough to show 
that 3. $\Longrightarrow$2. 

$\log q_{k+1}={\it o}\left( q_k \log q_k \right)$ means:
\begin{equation*}
\forall \epsilon >0 \; \; \exists \Hat{n} \left( \epsilon \right) \text{ such 
that } \forall l > \Hat{ n }\left( \epsilon \right) \quad \frac{\log 
q_{l+1}}{q_l \log q_l}<\epsilon
\end{equation*}

\indent
If $\log q_{l+1}<a q_l^\alpha$ for some positive constants 
$a$ and $\alpha < 1$ then:
\begin{equation*}
\frac{1}{\log q_k}\sum_{l=0}^k \frac{\log q_{l+1}}{q_l} \leq \frac{a}{\log 
q_k}\sum_{l=0}^\infty \frac{1}{q_l^{1-\alpha}} \leq \frac{aC}{\log q_k} 
\end{equation*}
for some universal constant $C$ thanks to (\ref{universalconstants}).

\indent
If $\log q_{l+1} \geq a q_l^\alpha$  and $\frac{1}{2}< 
\alpha <1$, consider the decomposition:
\begin{equation}
\frac{1}{\log q_k}\sum_{l=0}^k \frac{\log 
q_{l+1}}{q_l}=\underbrace{\frac{\log q_{k+1}}{q_k\log q_k}}_1 
+\underbrace{\frac{1}{\log q_k}\sum_{l=0}^{\Hat{n}\left( \epsilon \right)} 
\frac{\log q_{l+1}}{q_l}}_2+\underbrace{\frac{1}{\log 
q_k}\sum_{l=\Hat{n}\left( \epsilon \right)+1}^{k-1} \frac{\log 
q_{l+1}}{q_l}}_3
\label{threeterms}
\end{equation}
if $k-1 \geq \Hat{n}\left( \epsilon \right)+1$ otherwise the second and 
the third terms are replaced by $\frac{1}{\log 
q_k}\sum_{l=0}^{k-1}\frac{\log q_{l+1}}{q_l}$.
The third term can be bounded from above by:
\begin{equation*}
\frac{1}{\log q_k}\sum_{l=\Hat{n}\left( \epsilon \right)+1}^{k-1} 
\frac{\log q_{l+1}}{q_l} \leq \frac{\epsilon}{\log 
q_k}\sum_{l=\Hat{n}\left( \epsilon \right)}+1^{k-1} \log q_{l} \leq 
\epsilon \left( k-1-\Hat{n}\left( \epsilon \right) 
\right)\frac{\log q_{k-1}}{\log q_k}\; . 
\end{equation*}
Since  $\log q_j \leq \frac{2}{e}q_j^{\frac{1}{2}}$, from (\ref{growthqn}) and 
the hypothesis $\log q_{l+1} \geq a q_l^\alpha$ we obtain:
\begin{align*}
\frac{1}{\log q_k}\sum_{l=\Hat{n}\left( \epsilon \right)+1}^{k-1} 
\frac{\log q_{l+1}}{q_l} &\leq \left( k-1 \right) \frac{\epsilon}{a 
q_{k-1}^\alpha}\frac{2}{e}q_{k-1}^{\frac{1}{2}} \leq \\
 &\leq \frac{\epsilon 2}{ea}\left( k-1 
\right)e^{-\left( k-2 \right) \left( \alpha -\frac{1}{2} \right)\log G} 
\leq \epsilon C_1
\end{align*}
with $C_1 =\frac{2}{ea} \frac{e^{-1+\left( 
\alpha-\frac{1}{2} \right) \log G}}{\left( \alpha -\frac{1}{2}\right) \log 
G}$, $G=\frac{\sqrt{5}+1}{2}$.

\indent

The second term of (\ref{threeterms}) is bounded by
\begin{equation*}
\frac{1}{\log q_k}\sum_{l=0}^{\Hat{n}\left( \epsilon \right)} \frac{\log 
q_{l+1}}{q_l} \leq \frac{C_2}{\left( k-1 
\right) \log G -\log 2} \le \epsilon C_{2}
\end{equation*}
if $k > k\left( \epsilon \right)>\hat{n}(\epsilon )$, for some positive constant $C_2$.

\indent
Putting these estimates together we can bound (\ref{threeterms}) with:
\begin{equation*}
\frac{1}{\log q_k}\sum_{l=0}^k \frac{\log q_{l+1}}{q_l} \leq \epsilon + 
\epsilon C_1+\epsilon C_2
\end{equation*}
for all $\epsilon >0$ and for all $k > 
k\left( \epsilon \right)$, thus
$\sum_{l=0}^k\frac{\log_{l+1}}{q_l}={\it o}\left( \log q_k \right)$
\endproof

\subsection{Divergence of the modified linearization power series 
when the artihmetical conditions of Theorem \ref{theultradifferentiable} are not 
satisfied}

In Theorem 
\ref{theultradifferentiable} we proved that if $F \in 
z\mathbb{C}\left\{ z \right\}$ and $\omega$ verifies condition 
(\ref{BrjunoM})
then the linearization  $H\in z\mathbb{C}\left[ \left[ z \right] \right]_{\left( 
M_n \right)}$. The power series coefficients $h_{n}$ of 
$H$ are given by (\ref{recursive}). 

Let us define the sequence of strictly positive real numbers 
$(\Tilde{h}_n)_{n\ge 0}$ as follows:
\begin{equation}
\Tilde{h}_0=1\; , \;\;
\Tilde{h}_n =\frac{1}{|\lambda^n-1|} \sum_{m=2}^{n+1}|f_{m}|
\sum_{n_{1}+\ldots +n_{m}= n+1-m\, , \, 
n_{i}\ge 0 } \Tilde{h}_{n_{1}}\ldots \Tilde{h}_{n_{m}}\; . 
\label{htilde}
\end{equation}
Clearly 
$|h_{n}|\le \Tilde{h}_{n+1}$. Let 
$\Tilde{H} $ denote the formal power series associated to the sequence 
$(\Tilde{h}_n)_{n\ge 0}$
\begin{equation}
\Tilde{H}(z) = \sum_{m=1}^{\infty}\Tilde{h}_{n-1}z^n
\label{Htilde}
\end{equation}
Following closely \cite{Yoccoz}, Appendice 2, 
in this section we will prove that if condition (\ref{BrjunoM}) is 
violated then $\Tilde{H}$ doesn't belong to $z\mathbb{C}\left[ \left[ z 
\right] \right] _{\left( M_n \right)}$.

Note that since it is not restrictive to assume that $|f_{2}|\ge 1$
one has 
\begin{equation}
	\label{lowerh}
	\Tilde{h}_{n}> \sum_{k=0}^{n-1}\Tilde{h}_{k}\Tilde{h}_{n-1-k}
	\ge \Tilde{h}_{n-1}\; , 
\end{equation}
	\indent
thus the sequence $(\Tilde{h}_{n})_{n\ge 0}$ is strictly increasing. 

Let $\omega$ be an irrational  number which violates (\ref{BrjunoM}) 
and let $U=\{ q_j:q_{j+1} \geq \left( q_j +1 \right)^2 \}$ where $\left( 
q_j \right)_{j \geq 1}$ are the denominators of the convergents of $x$. 
Since $\inf_n \frac{1}{n}\log M_n=c>-\infty$  we have:
\begin{equation*}
\sum_{q_j \not\in U,j=0}^{k\left( n \right)}\frac{\log 
q_{j+1}}{q_j}-\frac{\log M_n}{n} \leq \sum_{q_j \not\in U,j=0}^{k\left( n 
\right)}\frac{2\log \left( q_j +1 \right)}{q_j}-c= \tilde{c} < +\infty
\end{equation*}
where $k\left( n \right)$ is defined by: $q_{k\left( n \right)} \leq n < 
q_{k\left( n \right)+1}$. 

On the other hand $\limsup_{n \rightarrow \infty} \left( 
\sum_{j=0}^{k\left( n \right)}\frac{\log q_{j+1}}{q_j} -\frac{\log M_n}{n} 
\right)= ~\infty$ thus
\begin{equation}
\limsup_{n \rightarrow \infty} \left( \sum_{q_j \in U:j=0}^{k\left( n 
\right)}\frac{\log q_{j+1}}{q_j} -\frac{\log M_n}{n} \right)= \infty
\end{equation}
this implies that $U$ is not empty. From now on the elements of $U$ will 
be denoted by: $q_0^\prime < q_1^\prime < \ldots$.

\indent
Let $n_i= \lfloor \frac{q_{i+1}^\prime}{q_{i}^\prime+1} \rfloor$. 

\begin{lemma}
The subsequence  $\left( \Tilde{h}_{q_{i}^\prime} \right)_{i \geq 0}$ 
verifies:
\begin{equation}
\label{subsequenceestimate}
\Tilde{h}_{q_{i+1}^\prime} \geq \frac{1}{ |\lambda^{q_{i+1}^\prime}-1| } 
\Tilde{h}_{q_{i}^\prime}^{n_i} \; . 
\end{equation}
\end{lemma}

\proof
From the definition (\ref{htilde}) and the assumption $|f_{2}|\ge 1$
it follows that
\begin{equation*}
	\Tilde{h}_{2s-1}\ge 
	\frac{|f_{2}|}{|\lambda^{2s-1}-1|}\Tilde{h}_{s-1}^{2}
	\ge \frac{\Tilde{h}_{s-1}^{2}}{2}
\end{equation*}
thus for all $i\ge 2$ and $s\ge 1$ one has
\begin{equation}
\label{geometric-h}
\Tilde{h}_{2s-1}\ge \frac{\Tilde{h}_{s-1}^{i}}{2}\; . 
\end{equation}
\indent
Choosing $s=q_{i}^\prime +1$, $i=n_{i}$ 
this leads to the desired 
estimate:
\begin{equation*}
	\Tilde{h}_{q_{i+1}^\prime}\ge \frac{2|f_{2}|}{|\lambda^{q_{i+1}^\prime}-1|}
	\Tilde{h}_{q_{i+1}^\prime-1}\ge 
	\frac{2|f_{2}|}{|\lambda^{q_{i+1}^\prime}-1|}
	\Tilde{h}_{n_{i}(q_{i}^\prime +1)-1}\ge 
	\frac{\Tilde{h}_{q_{i}^\prime}^{n_i}}{|\lambda^{q_{i+1}^\prime}-1|}\; . 
\end{equation*}
\endproof

\indent
By means of the previous lemma we can now prove that 
$\limsup_{n \rightarrow \infty}\frac{1}{n}\log 
\frac{\Tilde{h}_n}{M_n}=+\infty$.

Let $\alpha_{i}=n_{i}\frac{q_{i}^\prime}{q_{i+1}^\prime}$. Then 
$1 \geq \alpha_{i}\geq 
\left( 1-\frac{1}{q_i^\prime +1}\right)^2$, which assures that $\prod_{i 
\geq 0}\alpha_{i} =c$ for some finite constant $c$ (depending on $\omega$).
Then from (\ref{subsequenceestimate}) we get:
\begin{equation*}
\frac{1}{q_{i+1}^\prime}\log 
\frac{\Tilde{h}_{q_{i+1}^\prime}}{M_{q_{i+1}^\prime}}\geq c \left[ 
\sum_{j=1}^{i+1}-\frac{\log |\lambda^{q_{j}^\prime}-1|}{q_j^\prime} 
-\frac{1}{q^\prime_{i+1}}\log M_{q^\prime_{i+1}}\right]+c_4 
\end{equation*}
which diverges as  $i\rightarrow \infty$.

\appendix

\section{continued fractions and Brjuno's numbers}
\label{appfrazioni}

Here we summarize briefly some basic notions on continued fraction 
development and we define the Brjuno numbers. 

\indent
For a real number $\omega$, we note $\lfloor \omega \rfloor$ its integer part 
and $\{ \omega\}=\omega - \lfloor \omega \rfloor$ its fractional part.
We define the Gauss' continued fraction algorithm:
\begin{itemize}
\item $a_0=\lfloor \omega \rfloor$ and $\omega_0=\{ \omega\}$
\item for all $n \geq 1$: $a_n=\lfloor \frac{1}{\omega_{n-1}} \rfloor$ and 
$\omega_n=\{ \frac{1}{\omega_{n-1}} \}$
\end{itemize}
namely the following representation of $\omega$:
\begin{equation*}
\omega=a_0+\omega_0=a_0+\frac{1}{a_1+\omega_1}=\ldots
\end{equation*}
For short we use the notation $\omega=\left[ 
a_0,a_1,\ldots,a_n,\ldots\right]$.

\indent
It is well known that to every expression $\left[ 
a_0,a_1,\ldots,a_n, \ldots \right]$ there corresponds a unique 
irrational number. Let 
us define the sequences $\left( p_n \right)_{n\in \mathbb{N}}$ and 
$\left( q_n \right)_{n\in \mathbb{N}}$ as follows:
\begin{eqnarray*}
q_{-2}=1\text{, }q_{-1}=0\text{, }q_n=a_n q_{n-1}+q_{n-2}\\
p_{-2}=0\text{, }p_{-1}=1\text{, }p_n=a_n p_{n-1}+p_{n-2}
\end{eqnarray*}
It is easy to show that: $\frac{p_n}{q_n}=\left[ 
a_0,a_1,\ldots,a_n\right]$.

\indent
For any given $\omega \in \mathbb{R} \setminus \mathbb{Q}$ the sequence 
$\left( \frac{p_n}{q_n} \right)_{n\in \mathbb{N}}$ satisfies
\begin{equation}
	\label{growthqn}
	q_{n}\ge \left(\frac{\sqrt{5}+1}{ 2}\right)^{n-1}\; , \; \; n \ge 1
\end{equation}
thus 
\begin{equation}
\label{universalconstants}
\sum_{k\geq 0}\frac 1{q_k}\leq \frac{\sqrt{5}+5}{2}\;\;\;\text{ and }\;\;\;\sum_{k\geq 
0}\frac{\log q_k}{q_k}\leq \frac 1e\frac{2^{\frac 54}}{2^{\frac 
34}-1}\; , 
\end{equation}
and it has the following important properties:
\begin{lemma}
for all $n \geq 1$ then: $\frac{1}{q_n+q_{n+1}} \leq \lvert q_n \omega - 
p_n \rvert < \frac{1}{q_{n+1}}$.
\end{lemma}
\begin{lemma}
\label{lemmanumb}
If for some integer $r$ and $s$, $\mid \omega -\frac {r}{s} \mid \leq 
\frac {1}{2 s^{2}}$ then $\frac {r}{s} = \frac {p_k}{q_k}$ for some 
integer $k$.
\end{lemma}
\begin{lemma}
\label{bestapproximation}
The law of best approximation: if $1\le q\le 
q_{k}$, $(p,q)\not= (p_{n},q_{n})$ and $n\ge 1$ then 
$|qx-p|>|q_{n}x-p_{n}|$. Moreover if $(p,q)\not= (p_{n-1},q_{n-1})$
then $|qx-p|>|q_{n-1}x-p_{n-1}|$.
\end{lemma}
For a proof of these standard lemmas we refer to 
\cite{Hardy}.

\indent
The growth rate of $\left( q_n \right)_{n\in \mathbb{N}}$ describes 
how rapidly $\omega$ can be approximated by rational numbers.
For example $\omega$ is a diophantine number \cite{Siegel} if and only 
if there 
exist two constants $c >0$ and $\tau \ge 1$ such that 
$q_{n+1}\le c q_{n}^\tau$ for all $n\ge 0$. 

\indent
To every $\omega \in \mathbb{R} \setminus \mathbb{Q}$ we associate, using 
its convergents, an arithmetical  function: 
\begin{equation}
	\label{Brjunofunction}
	B\left( \omega \right)= \sum_{n \geq 0} 
\frac{\log q_{n+1}}{q_n}
\end{equation}
We say that $\omega$ is a {\it Brjuno number} or 
that it satisfies the {\it Brjuno condition} if $B\left( \omega \right)<+\infty$. 
The Brjuno  condition gives a limitation to the growth rate of 
$\left( q_n \right)_{n\in \mathbb{N}}$. It was originally introduced by 
A.D.Brjuno \cite{Brjuno}. 
The  Brjuno condition is weaker than the 
Diophantine condition: for example if $a_{n+1} \le c e^{a_{n}}$ for 
some positive constant $c$ and for all $n\ge 0$ then $\omega = 
[a_{0}, a_{1}, \ldots , a_{n}, \ldots ]$ is a Brjuno number but is not 
a diophantine number.

\section{Davie's lemma}
\label{appDavie}

In this appendix we summarize the result of \cite{Davie} that we use, 
in particular  Lemma \ref{lemmaDavie}.
Let $\omega \in \mathbb{R} \setminus \mathbb{Q}$ and $\left\{ q_n 
\right\}_{n \in \mathbb{N}}$ the partial denominators of the 
continued fraction for $\omega$ in the Gauss' development. 

\begin{definition}
Let $A_k = \left\{ n \geq 0 \mid \| n \omega \| \leq \frac{1}{8q_k} 
\right\}$, $E_k=\max \left( q_k , q_{k+1}/4 \right)$ and $\eta_k = 
q_k / E_k$.
Let $A_k^{*}$ be the set of non negative integers 
$j$ such that either $j \in A_k$ or for some $j_1$ and $j_2$ in 
$A_k$, with $j_2-j_1 < E_k$, one has $j_1 < j < j_2$ and $q_k$ 
divides $j-j_1$.
For any non negative integer $n$ define:
\begin{equation*}
l \left( n \right) = \max \left\{ \left( 1+\eta_k \right) 
\frac{n}{q_k}-2 , \left( m_n \eta_k+n \right) \frac{1}{q_k}-1 \right\}
\end{equation*}
where $m_n = \max \{ j \mid 0 \leq j \leq n , j \in A_k^{*} \}$.
We then define the function $h_k \left( n \right)$
\begin{equation*}
h_k \left( n \right)= \begin{cases}
				\frac{m_n+\eta_k n}{q_k}-1& \text{if $m_n+q_k \in A_k^{*}$} \\
				l \left( n \right)& \text{if $m_n+q_k \not\in A_k^{*}$}
			\end{cases}
\end{equation*}
\end{definition}
The function $h_k \left( n \right)$ has some properties collected in 
the following proposition

\begin{proposition}
The function $h_k \left( n \right)$ verifies;
\begin{enumerate}
\item $\frac{\left( 1+\eta_k \right)n}{q_k}-2 \leq h_k \left( n 
\right) \leq \frac{\left( 1+\eta_k \right)n}{q_k}-1$ for all $n$.
\item If $n>0$ and $n \in A_k^{*}$ then $h_k \left( n \right) \geq 
h_k \left( n -1 \right)+1$.
\item $h_k \left( n \right) \geq h_k \left( n-1 \right)$ for all 
$n>0$.
\item $h_k \left( n+q_k \right) \geq h_k \left( n \right) +1$ for all 
$n$.
\end{enumerate}
\end{proposition}

\indent
Now we set $g_k \left( n \right)= \max \left( h_k \left( n \right), 
\lfloor \frac{n}{q_k} \rfloor \right)$ and we state the following 
proposition

\begin{proposition}
\label{gpropos}
The function $g_k $ 
is non negative and verifies:
\begin{enumerate}
\item $g_k \left( 0 \right)=0$
\item $g_k \left( n \right) \leq \frac{\left( 1+\eta_k 
\right)n}{q_k}$ for all $n$
\item $g_k \left( n_1 \right) + g_k \left( n_2 \right) \leq g_k 
\left( n_1 +n_2 \right)$ for all $n_1$ and $n_2$
\item if $n \in A_k$ and $n>0$ then $g_k \left( n \right) \geq g_k 
\left( n-1 \right)+1$
\end{enumerate}
\end{proposition}
The proof of these propositions can be found in \cite{Davie}. 

Let $k(n)$ be defined by the condition $q_{k(n)}\le n <q_{k(n)+1}$. 
Note that $k$ is non--decreasing. 

\begin{lemma} {\bf Davie's lemma}
	\label{lemmaDavie}
	Let 
\begin{equation*}
	K(n)=n\log 2+\sum_{k=0}^{k(n)}g_{k}(n)\log(2q_{k+1})\; . 
\end{equation*}
The function $K \left( n \right)$ verifies:
\begin{enumerate}
\item There exists a universal constant $\gamma_{3}>0$ such that 
\begin{equation*}
K(n)\le n\left(\sum_{k=0}^{k(n)}\frac{\log 
q_{k+1}}{q_{k}}+\gamma_{3}\right)\; ; 
\end{equation*}
\item $K(n_{1})+K(n_{2})\le K(n_{1}+n_{2})$ for all $n_{1}$ and $n_{2}$; 
\item $-\log |\lambda^n -1| \le K(n)-K(n-1)$. 
\end{enumerate}
\end{lemma}

The proof is a straightforward application of Proposition \ref{gpropos}.

\end{document}